\documentclass[11pt]{amsart}
\usepackage{amssymb,amsthm}

\theoremstyle{plain}
\newtheorem{theorem}{Theorem}
\newtheorem{prop}{Proposition}
\theoremstyle{remark}
\newtheorem{rem}{Remark}
\newtheorem{pro}{Problem}

\newcommand{\N}{\mathbb{N}}
\newcommand{\cD}{\mathcal{D}}
\newcommand{\cR}{\mathcal{R}}

\newcommand{\cU}{\mathcal{U}}

\newcommand{\eps}{\varepsilon}
\newcommand{\Diff}{\mathrm{Diff}_\mu^{1+\alpha}(M)}
\newcommand{\Diffum}{\mathrm{Diff}_\mu^1(M)}
\newcommand{\SE}{\mathcal{SE}}
\begin{document}

\title{A remark on conservative diffeomorphisms}
%\date{\today}
\author[J.~Bochi, B.~R.~Fayad, E.~Pujals]{Jairo Bochi, Bassam R.~Fayad, and Enrique Pujals}
\thanks{J.B.~was supported by CNPq-Profix during the preparation of this work. 
J.B.~thanks the hospitality of the LAGA -- Universit\'e de Paris 13.}

\address{Inst.~Matem\'atica -- UFRGS -- Av Bento Gon\c{c}alves 9500 -- 91509-900 Porto Alegre -- Brazil.}
\email{jairo@mat.ufrgs.br}

\address{LAGA -- Universit\'e Paris 13 -- 99 Av J-B. Cl\'ement -- 93430 Villetaneuse -- France.}
\email{fayadb@math.univ-paris13.fr}

\address{IMPA -- Estr.~D.~Castorina 110 -- 22460-320 Rio de Janeiro -- Brazil.}
\email{enrique@impa.br}

\begin{abstract}
We show that a stably ergodic diffeomorphism can be $C^1$ approximated by
a diffeomorphism having stably non-zero Lyapunov exponents.

\vskip 1 \baselineskip
{\bf Une remarque sur les diff\'eomorphismes conservatifs}

\vskip 0.5\baselineskip
\noindent {\sc R\'esum\'e.} On montre qu'un diff\'eomorphisme stablement 
ergodique peut \^etre $C^1$ approch\'e par un diff\'eomorphisme ayant des
exposants de Lyapunov stablement non-nuls.

\end{abstract}

\maketitle

%%%%%%%%%%%%%%%%%%%%%%%%%%%%%%%%%%%%%%%%%%%%%%%

Two central notions in Dynamical Systems are ergodicity and hyperbolicity.
In many works showing that certain systems are ergodic,
some kind of hyperbolicity (e.g.~uniform, non-uniform or partial)
is a main ingredient in the proof.
In this note the converse direction is investigated.

\medskip

Let $M$ be a compact manifold of dimension $d \ge 2$, and let $\mu$ be a volume measure in $M$.
Take $\alpha>0$ and let $\Diff$ be the set of $\mu$-preserving $C^{1+\alpha}$ diffeomorphisms,
endowed with the $C^1$ topology.
Let $\SE \subset \Diff$ be the set of stably ergodic diffeomorphisms
(i.e., the set of diffeomorphisms such that every sufficiently
$C^1$-close $C^{1+\alpha}$ conservative diffeomorphism is ergodic).

Our result answers positively a question of~\cite{BuDoPe}:

\begin{theorem}\label{t.1}
There is an open and dense set $\cR \subset \SE$ such that if $f \in \cR$ then
$f$ is non-uniformly hyperbolic, that is, all Lyapunov exponents of $f$ are non-zero.
Moreover, every $f\in \cR$ admits a dominated splitting
$TM = E^+ \oplus E^-$, where $E^+$ (resp.~$E^-$)
coincides a.e.~with the sum of the Oseledets spaces corresponding to positive (resp.~negative)
Lyapunov exponents.
\end{theorem}

\begin{rem}
The set $\SE$ contains all Anosov diffeomorphisms, and many partially hyperbolic ones -- see e.g.~\cite{GPS}.
It is not true that every stably ergodic diffeomorphism can be approximated by a
partially hyperbolic system, see ~\cite{T, BnV}.
\end{rem}

\begin{rem}
Let $\SE'$ be the set of diffeomorphisms $f\in\SE$
such that every power $f^k$, $k \ge 2$, is ergodic.
Then every $f$ in $\SE' \cap \cR$ is Bernoulli.
This follows from theorem~\ref{t.1} and Pesin theory
(see theorem 5.10 in \cite{L}).
\end{rem}

The proof of theorem~\ref{t.1} has three steps:
\begin{itemize}
\item[1.]
A stably ergodic (or stably transitive) diffeomorphism $f$ must have a dominated splitting.
This is true because if it doesn't, \cite{BDP}
permits us to perturb $f$ and create a periodic point whose derivative is the identity.
Then, using the Pasting Lemma from~\cite{AM}
(for which $C^{1+\alpha}$ regularity is an essential hypothesis),
one breaks transitivity.

\item[2.]
A result of~\cite{BB} gives a perturbation of $f$ such that the sum of the Lyapunov exponents
``inside'' each of the bundles of the (finest) dominated splitting is non-zero.

\item[3.]
Using a result of~\cite{BV}, we find another perturbation
such that the Lyapunov exponents in each of the bundles become almost equal.
(If we attempted to make the exponents exactly equal, we couldn't guarantee that the
perturbation is $C^{1+\alpha}$.)
Since the sum of the exponents in each bundle varies continuously,
we conclude there are no zero exponents.
\end{itemize}

\begin{rem}
The perturbation techniques of \cite{BB} and \cite{BV} in fact don't assume
ergodicity, but are only able to control the integrated Lyapunov exponents.
That's why we have to assume stable ergodicity (in place of stable transitivity)
in theorem~\ref{t.1}.
\end{rem}

\begin{rem}
Theorem~\ref{t.1} is stated in $C^1$ topology because in higher topologies
the technology from~\cite{BDP}, \cite{BB}, and \cite{BV} is not available.
The $C^{1+\alpha}$ diffeomorphisms come from~\cite{AM}.
To get our result in $C^1$ topology (which perhaps would be more natural)
one has to solve the following problem: 
any diffeomorphism having a periodic point tangent to the identity may be
$C^1$-approximated by a non-transitive diffeomorphism.
\end{rem}

\begin{rem}
Some ideas of the present proof were already present in~\cite{DP}.
\end{rem}

%%%%%%%%%%%%%%%%%%%%%%%%%%%%%%%%%%%%%%%%%%%%
\medskip

Let us recall briefly the definition and some properties of dominated splittings,
see~\cite{BDP} for details.
Let $f\in \Diffum$.

A $Df$-invariant splitting $T M = E^1 \oplus \cdots \oplus E^k$, with $k\ge 2$, is
called a \emph{dominated splitting} (over $M$)
if there are constants $c$, $\tau>0$ such that
\begin{equation}\label{e.def ds}
\frac{\|Df^n(x)\cdot v_j\|}{\|Df^n(x)\cdot v_i\|} < c e^{-\tau n}
\end{equation}
for all $x\in M$, all $n\ge 1$, and all unit vectors
$v_i\in E^i(x)$ and $v_j\in E^j(x)$, provided $i<j$.
(One can also define in the same way a dominated splitting over an $f$-invariant set.)

A dominated splitting is always continuous, that is, 
the spaces $E_i(x)$ depend continuously on $x$.
Also, a dominated splitting persists under $C^1$-perturbations of the map.
More precisely, if $g$ is sufficiently close to $f$, then $g$ has 
a dominated splitting $E^1_g \oplus \cdots \oplus E^k_g$,
called the \emph{continuation},
with $\dim E^i_g = \dim E^i$ and which coincides with the given one when $g=f$.
Moreover, $E^i_g(x)$ depends continuously on $g$ (and $x$).

A dominated splitting $E^1 \oplus \cdots \oplus E^k$ is called
the \emph{finest dominated splitting} if there is no dominated splitting defined over
all $M$ with more than $k$ bundles.
If some dominated splitting exists, then the finest dominated splitting exists,
is unique, and refines every dominated splitting.

The continuation of the finest dominated splitting
is not necessarily the finest dominated splitting of the perturbed diffeomorphism.
We call a dominated splitting for $f \in \Diff$ \emph{stably finest} if 
it has a continuation which is the finest dominated splitting 
of every sufficiently $C^1$-close diffeomorphism of class $C^{1+\alpha}$.
It is easy to see that diffeomorphisms with stably finest dominated splittings are
(open and) dense among $C^{1+\alpha}$ diffeomorphisms with a dominated splitting.

%%%%%%%%%%%%%%%%%%%%%%%%%%%%%%%%%%%%%%
\medskip

Let $\lambda_1(f,x)\ge\cdots\ge\lambda_d(f,x)$ be the Lyapunov exponents of $f$ (counted with multiplicity),
defined for almost all~$x$.
(See e.g.~\cite{A} for definition and basic properties of Lyapunov exponents.)
We write also
\begin{equation}\label{e.ie}
\lambda_i(f) = \int_M \lambda_i(f,x)\,d\mu(x).
\end{equation}

Assume $f$ has a dominated splitting $E^1 \oplus \cdots \oplus E^k$.
Then the Oseledets splitting is a measurable refinement of it.
For simplicity of writing, we will say the exponent $\lambda_p$ \emph{belongs} to the bundle $E^i$
if
$d_1 + \cdots + d_{i-1} < p \le d_1 + \cdots + d_i$,
where $d_i = \dim E^i$.
By~\eqref{e.def ds}, there is an uniform gap between Lyapunov exponents that belong to different bundles.

%%%%%%%%%%%%%%%%%%%%%%%%%%%%%%%%%%%%%%%%%%%%
\medskip

We now give the proof of theorem~\ref{t.1} in detail.
Let $\cR$ be the set of $f \in \SE$ such that
$f$ has a dominated splitting $E^+ \oplus E^-$
with $\lambda_p(f)>0>\lambda_{p+1}(f)$, where $p=\dim E^+$.
First we see that $\cR$ is an open set.
Indeed, given $f\in \cR$, there is an open set $\cU \ni f$
where the dominated splitting has a continuation,
say $E^+_g \oplus E^-_g$ for $g\in \cU$.
As $\lambda_{p+1}$ is the top exponent in $E^-$, we can write
\begin{equation}\label{e.top}
\lambda_{p+1} (g) = \inf_{n \in \N} 
\frac{1}{n} \int_M \log \|Dg^n(x)|_{E_g^-(x)} \| \; d\mu(x).
\end{equation}
Therefore $g \in \cU \mapsto \lambda_{p+1}(g)$ is an upper semicontinuous function.
Accordingly, $\lambda_{p+1}(g) < 0$ for all $g$ sufficiently close to $f$.
And analogously for $\lambda_p$, showing that $\cR$ is open.

\medskip

Next we show that $\cR$ is dense in $\SE$.
Take $f\in\Diff$ a stably ergodic diffeomorphism.
As mentioned, this implies that $f$ has a dominated splitting, see~\cite{AM}.
As remarked above, we can assume, after a perturbation of $f$ if necessary,
that $f$ has a stably finest dominated splitting.

For all $g$ sufficiently close to $f$, we denote by  
$E^1_g \oplus \cdots \oplus E^k_g$ the finest dominated splitting of $g$.
Let us indicate by $J_i(g)$ the sum of all Lyapunov exponents $\lambda_p(g)$ that belong to $E^i_g$.
Then we can also write
\begin{equation} \label{e.j}
J_i (g) =  \int_M \log \big| \det Dg|_{E^i_g} \big| \, d\mu.
\end{equation}
In particular, $J_i(\mathord{\cdot})$ is a continuous function in the neighborhood of $f$.

By the theorem from~\cite{BB}, up to $C^1$-perturbing $f$,
we may assume $J_i(f)\neq 0$ for all $i$.
(It is important to notice that the perturbed map can be taken 
of class $C^{1+\alpha}$ since so is the original $f$.)

In the last step we need the following proposition:

\begin{prop}\label{p.almost equal}
Let $f\in \SE$.
Assume that $f$ has a stably finest dominated splitting $E^1_f\oplus\cdots\oplus E^k_f$.
Then for all $\eps>0$ there exists a perturbation $g\in\Diff$ of $f$ such that
if the Lyapunov exponents $\lambda_p(g)$, $\lambda_q(g)$ belong to the same bundle $E^i_g$,
then $|\lambda_p(g)-\lambda_q(g)|<\eps$.
\end{prop}

Applying the proposition, we find $g$ close to $f$ such that
all $\lambda_p(g)$ in $E^i_g$ are close to $J_i(g)/\dim E^i$ and therefore are non-zero.
This finishes the proof of theorem~\ref{t.1}, modulo giving the:

\begin{proof}[Proof of proposition~\ref{p.almost equal}]
For $f\in\Diff$ and $1\le p \le d$, let us write $\Lambda_p(f) = \lambda_1(f) + \cdots + \lambda_p(f)$.
Then $\Lambda_p(\mathord{\cdot})$ is an upper semicontinuous function 
(see~\cite{A} or \cite{BV}).
Since $\Diff$ is not a complete metric space, 
we can't deduce that the set of continuity points of
$\Lambda_p(\mathord{\cdot})$ is dense.
Nevertheless, for every $\eps>0$, the set
$$
\cD_{\eps,p} =
\{f \in \Diff;\; \exists \; \cU \ni f \text{ open s.t. }
|\Lambda_p(g_1) - \Lambda_p(g_2)| < \eps \; \forall g_1, g_2\in \cU\}
$$
is (open and) dense in $\Diff$.
(This is an easy exercise using $\Lambda_p \ge 0$.)
In particular, $\cD_\eps = \bigcap_{p=1}^d \cD_{\eps,p}$ is dense.

Now let $f\in \SE$ have a stably finest dominated spitting into $k$ bundles.
Fix $\eps>0$ and take $g \in \cD_\eps$ very $C^1$-close to $f$.
We claim that $g$ has the desired properties:
for any $i=1,\ldots,k$, if $\lambda_p$, $\lambda_q$ belong to $E^i_g$
then $\lambda_p$, $\lambda_q$ are close.
Clearly, it suffices to consider the case $q=p+1$.

Consider the set $D_p(g)$ of points $x\in M$ 
such that there exists a dominated splitting 
$T_{\overline{o(g,x)}} M = F \oplus G$ 
over the closure of the $g$-orbit of $x$, with $\dim F = p$.
Notice there is no dominated splitting 
$TM = F \oplus G$ (over $M$) with $\dim F = p$,
because $\lambda_p$ and $\lambda_{p+1}$ belong to the same bundle 
of the finest dominated splitting of $g$.
Thus no $x \in D_p(g)$ can have a dense orbit.
In particular, $D_p(g)$ has zero measure.
By proposition~4.17 from~\cite{BV}, there exists a
$C^1$-perturbation $h$ of $g$ such that
\begin{align*}
\Lambda_p(h) &< \Lambda_p(g) 
               - \int_{M \setminus D_p(g)} 
                    \frac{\lambda_p(g,x) - \lambda_{p+1}(g,x)}{2} \; d\mu(x)
               + \eps \\
               & = \Lambda_p(g) - \frac{\lambda_p(g) - \lambda_{p+1}(g)}{2} + \eps.
\end{align*}
(In the notation of~\cite{BV}, $\Gamma_p(g,\infty) = M \setminus D_p(g)$.)
Because $g$ is $C^{1+\alpha}$, the map $h$
given by the proof of proposition~4.17 in~\cite{BV} is $C^{1+\alpha}$ as well.
Since $g\in \cD_{\eps,p}$ and $h$ is close to $g$,
we have $|\Lambda_p(h) - \Lambda_p(g)| < \eps$ and accordingly
$\lambda_p(g) - \lambda_{p+1}(g) < 4\eps$.
\end{proof}

%%%%%%%%%%%%%%%%%%%%%%%%%%%%%%%%%%%%%%%%
\medskip

We close this note with some questions about what can be said
in the absence of stable ergodicity.
The following question (similar to one in \cite{SW})
is likely to have a positive answer:

\begin{pro}\label{p.1}
Is it true that for the generic $f\in\Diffum$, either all Lyapunov exponents
are zero at almost every point, or $f$ is non-uniformly hyperbolic (i.e.,
all Lyapunov exponents are non-zero almost everywhere)?
\end{pro}

Notice this is true if $\dim M=2$, by~\cite{B} (later extended in \cite{BV}).
Using the main result of the papers~\cite{BV} and~\cite{BB}, it is not difficult
to show that the dichotomy of problem~\ref{p.1} holds true
modulo an eventual positive answer to the following well known conjecture of A.~Katok:
\begin{pro}\label{p.2}
Is it true that the generic map $f\in\Diffum$ is ergodic?
\end{pro}

%\begin{rem}
%The theorem of Oxtoby-Ulam~\cite{OU} says that $C^0$-generic volume-preserving
%homeomorphisms are ergodic.
%Also, it was recently shown by Bonatti and Crovisier~\cite{BC}
%that the generic $f\in\Diffum$ is transitive.
%\end{rem}

%%%%%%%%%%%%%%%%%%%%%%%%%%%%%%%%%%%%%%%%%%%%%%%

\end{document}